\documentclass[12pt]{article}
\usepackage{a4wide,color,hyperref}
\usepackage{amssymb}
\usepackage{amsmath}

\newtheorem{theorem}{Theorem}
\newtheorem{question}{Question}

\newtheorem{corollary}{Corollary}
\newtheorem{lemma}{Lemma}

\renewcommand{\Bbb}[1]{\mathbb{#1}}
\newcommand{\N}{{\Bbb N}}         
\newcommand{\R}{{\Bbb R}}         
\newcommand{\Rp}{[0,+\infty)}    
\newcommand{\Z}{{\Bbb Z}}         

\newcommand{\Anm}{\cA_{n,m}}

\newcommand{\AnmP}{\cA_{n,m}(\Psi)}

\newcommand{\AdnmP}{\cA'_{n,m}(\Psi)}

\newcommand{\AddnmP}{\cA''_{n,m}(\Psi)}
\newcommand{\Ap}{\cA(\psi)}

\newcommand{\Adp}{\cA'(\psi)}

\newcommand{\Fnm}{\cF_{n,m}}

\newcommand{\FnmP}{\cF_{n,m}(\Psi)}

\newcommand{\FdnmP}{\cF'_{n,m}(\Psi)}

\newcommand{\FddnmP}{\cF''_{n,m}(\Psi)}

\newcommand{\Wnm}{\cW_{n,m}}

\newcommand{\Bnm}{\cB_{n,m}}

\newcommand{\BnmP}{\cB_{n,m}(\Psi)}

\newcommand{\anyA}{\cA^{\circ}_{n,m}(\Psi)}
\newcommand{\anyF}{\cF^{\circ}_{n,m}(\Psi)}
\newcommand{\anyW}{\cW^{\circ}_{n,m}(\Psi)}
\newcommand{\anyB}{\cB^{\circ}_{n,m}(\Psi)}

\def\KG{Khintchine-Groshev}

\newcommand{\cA}{{\cal A}}
\newcommand{\cB}{{\cal B}}
\newcommand{\cC}{{\cal C}}
\newcommand{\cD}{{\cal D}}

\newcommand{\cF}{{\cal F}}
\newcommand{\cG}{{\cal G}}

\newcommand{\cS}{{\cal S}}

\newcommand{\cU}{{\cal U}}

\newcommand{\cW}{{\cal W}}

\newcommand{\ve}{\varepsilon}

\newcommand{\diam}{\operatorname{diam}}
\newcommand{\dist}{\operatorname{dist}}

\newcommand{\vv}[1]{{\mathbf{#1}}}

\renewcommand{\le}{\leqslant}
\renewcommand{\ge}{\geqslant}
\newcommand{\nz}{\smallsetminus\{0\}}
\newcommand{\bnz}{\smallsetminus\{{\vv0}\}}

\begin{document}

\title{A note on zero-one laws in metrical \\ Diophantine approximation}

\author{Victor Beresnevich\footnote{EPSRC Advanced Research Fellow, grant EP/C54076X/1}
\\ {\small\sc York} \and Sanju Velani\footnote{Research supported by EPSRC  grant EP/E061613/1 and
INTAS grant 03-51-5070}
\\ {\small\sc York}}

\date{\small\it Dedicated to Wolfgang Schmidt on the occasion of his 75th birthday}

\maketitle



\vspace*{-4ex}

\section{Introduction}

Given $\psi:\N\to\Rp$, let $\cA(\psi)$ denote the set of
$x\in[0,1]$ such that
\begin{equation}\label{e:001+}
 |qx+p|<\psi(q)
\end{equation}
holds for infinitely many $(p,q)\in\Z\times\Z\nz$. In 1924,
Khintchine \cite{Khintchine-1924} established a beautiful and
strikingly simply criterion for the `size' of  $\cA(\psi)$ expressed
in terms of Lebesgue measure. Under the condition that $\psi$ is
monotonic, Khintchine's theorem states that the measure of
$\cA(\psi)$  is one (respectively, zero) if the sum $\sum_q \psi(q)
$ diverges (respectively, converges).  The monotonicity condition is
only required in the divergence  case  and moreover it is absolutely
crucial.  Duffin and Schaeffer \cite{Duffin-Schaeffer-41:MR0004859}
constructed a non-monotonic function $\psi$ for which $\sum_q
\psi(q) $ diverges but $\cA(\psi)$  is of zero  measure.  In other
words, without the monotonicity assumption, Khintchine's theorem is
false and the famous Duffin-Schaeffer conjecture provides the
appropriate statement.  The key difference is that in
(\ref{e:001+}), we impose coprimality on the integers $p$ and $q$.
Let $\Adp$ denote the resulting subset of $\cA(\psi)$.  The
Duffin-Schaeffer conjecture states that the measure of $\Adp$ is one
(respectively, zero) if the sum $\sum_q \varphi(r) \, \psi(q) \,
q^{-1} $ diverges (respectively, converges).
Although various partial results have been obtained, the full
conjecture represents a key unsolved problem in metric number
theory -- see
\cite{Beresnevich-Bernik-Dodson-Velani-Roth,Harman-1998a} for
details. Returning to the raw set $\cA(\psi)$, without
monotonicity and coprimality the appropriate analogue of
Khintchine's theorem has been formulated by Catlin
\cite{Catlin-76:MR0417098}. The Catlin conjecture  also remains
open. 

The upshot of the above discussion is that currently we are unable
to prove analogues of Khintchine's theorem for either of the
fundamental sets $\Ap$ and $\Adp$.  However , it is known that the
Lebesgue measure of  $\Ap$ and $\Adp$ is either 0 or 1. In the
case of $\Ap$ this zero-one law is due to Cassels
\cite{Cassels-50:MR0036787} and in the case of $\Adp$ it is due to
Gallagher \cite{Gallagher-61:MR0133297}. The goal of this note is
to establish the higher dimensional analogues of these classical
zero-one laws. For a discussion concerning the higher dimensional
analogues of the conjectures of Duffin-Schaeffer and Catlin see
\cite{Beresnevich-Bernik-Dodson-Velani-Roth}.

Throughout, $m \geq1 $ and $n \geq 1 $ are integers. Given
$\Psi:\Z^m\to\Rp$, let $\AnmP$ be the set of $\vv X\in[0,1]^{nm}$
such that
\begin{equation}\label{e:001}
 |\vv q\vv X+\vv p|<\Psi(\vv q)
\end{equation}
holds for infinitely many $(\vv p,\vv q)\in\Z^m\times\Z^n\bnz$.
Here $|\cdot|$ denotes the supremum norm in $\R^m$, $\vv X$ is
regarded as an $n\times m$ matrix and $\vv q$ is regarded as a
row. Thus, $\vv q\vv X\in\R^m$ represents a system of $m$ real
linear forms in $n$ variables. In higher dimensions the set $\Adp$
has two natural generalizations:
$$
\AdnmP=\{\vv X\in[0,1]^{nm}:\text{(\ref{e:001}) holds for i.m. }(\vv
p,\vv q)\text{ with }\gcd(\vv p,\vv q)=1\}\,,
$$
$$
\AddnmP=\{\vv X\in[0,1]^{nm}:\text{(\ref{e:001}) holds for i.m.
}(\vv p,\vv q)\text{ with }\gcd(p_i,\vv q)=1\ \
j=\overline{1,m}\}\,.
$$

\noindent Here `i.m.' stands for `infinitely many' and   $\gcd(\vv
p,\vv q)$ denotes the greatest common divisor of all the
components of $\vv p$ and $\vv q$. If $\gcd(\vv p,\vv q)=1$ then
we say that $\vv p$ and $\vv q$ are \emph{coprime}.

Before we state our main result, let us agree on the following
notation:  $\anyA$ will denote any of the fundamental  sets $\AnmP$,
$\AdnmP$ and $\AddnmP$. Thus, a statement  for $\anyA$ is valid for
$\AnmP$, $\AdnmP$ and $\AddnmP$. Also, $|X|$ will denote the
$k$-dimensional Lebesgue measure of the set $X\subset\R^{k} $.

\begin{theorem}\label{t1}
For any $n,m$ and $\Psi$ we have that
$|\anyA|\in\{0,1\}$\,.
\end{theorem}

\section{Auxiliary results}

In this section we group together various self contained
statements  that we appeal to during the course of establishing
Theorem 1. Most are higher dimensional analogues of well known
one-dimensional statements.    Indeed, the one-dimensional version
of our  first result can be found in  \cite{Cassels-50:MR0036787}.

\begin{lemma}\label{lemma0}
Let $\{B_i\}$ be a sequence of balls in $\R^k$ with $|B_i|\to0$ as
$i\to\infty$. Let $\{U_i\}$ be a sequence of Lebesgue measurable
sets such that $U_i\subset B_i$ for all $i$. Assume that for some
$c>0$,
\begin{equation}\label{e:005}
     |U_i|\ge c|B_i|\qquad\text{for all }i\,.
\end{equation}
Then the sets
$$
\textstyle
 \cU=\limsup\limits_{i\to\infty}U_i:=\bigcap\limits_{j=1}^\infty\ \bigcup\limits_{i\ge
 j}U_i\qquad\text{ and }\qquad
 \cB=\limsup\limits_{i\to\infty}B_i:=\bigcap\limits_{j=1}^\infty\ \bigcup\limits_{i\ge
 j}B_i
$$
have the same Lebesgue measure.
\end{lemma}

\noindent\textit{Proof.} Let $\cU_j:=\bigcup_{i\ge j}U_i$ and
$\cD_j:=\cB\setminus\cU_j$. Then,
$\cD:=\cB\setminus\cU=\bigcup_{j}\cD_j$ and Lemma~\ref{lemma0}
states that $\cD$ has measure zero. Equivalently, every $\cD_j$
must have zero measure. Assume the contrary. Then, there is an
$l\in\N$ such that $|D_l|>0$ and therefore there is a density
point ${\vv x}_0$  of $\cD_l$. Since ${\vv x}_0\in\cB$, we have
that $x_0 \in B_{j_i}$ for a sequence $j_i$. Since
$|B_{j_i}|\to0$, we have that $|\cD_l\cap B_{j_i}|\sim|B_{j_i}|$
as $i\to\infty$. Since $\cD_j\supset \cD_l$ for all $j\ge l$,  it
follows that
\begin{equation}\label{e:006}
    |\cD_{j_i}\cap B_{j_i}|\sim|B_{j_i}|\text{ as }i\to\infty\,.
\end{equation}
On the other hand, by construction $\cD_{j_i}\cap
U_{j_i}=\emptyset$.  Thus, in view of (\ref{e:005}) we have that
$$
|B_{j_i}|\ge|U_{j_i}|+|\cD_{j_i}\cap B_{j_i}|\ge
c|B_{j_i}|+|\cD_{j_i}\cap B_{j_i}| \ ;
$$
i.e. $ |\cD_{j_i}\cap B_{j_i}|\le(1-c) \ |B_{j_i}| $ for all
sufficiently large $i$.  This contradicts (\ref{e:006}).

\vspace{-3ex}

\hfill $\boxtimes$

\vspace{1.5ex}

The following lemma is the higher dimensional analogue of the well
know one-dimensional `ergodic' property of rational
transformations -- see for example
\cite[Lemma~3]{Gallagher-61:MR0133297},
\cite[Lemma~2.2]{Harman-1998a} or \cite[Lemma~7]{Sprindzuk-1979-Metrical-theory}.

\begin{lemma}\label{l4}
 For any integer $l\ge2$ and $\vv s\in\Z^k$ consider the
transformation of the unit cube $ [0,1]^k $ into itself  given by
\begin{equation*}
   T \, : \,  \vv x\mapsto l \, \vv x+\frac{1}{l}\,\vv s\pmod 1 \ .
\end{equation*}
Let $A$ be a subset of $[0,1]^k$ such that $T(A) \subseteq A$.
Then $A$ is of Lebesgue measure $0$ or $1$.
\end{lemma}

\noindent\textit{Proof.} Let $A$ be as in the statement.   Then
$T^{\nu}(A) \subseteq A$,   where $
 T^{\nu}:    \vv x\mapsto l^\nu\vv x+\frac{\vv s}{l}\pmod 1
$ is the $\nu$-th iterate  of $T$.  Let $\chi_A$ be the
characteristic function of $A$. It follows that
\begin{equation}\label{e:008}
    \chi_A(\vv x)\le\chi_A\Big(l^\nu\vv x+\frac{\vv s}{l}\Big)\,.
\end{equation}
Suppose that $|A|>0$. Then there is a density point $\vv x_0$ of
$A$. Let $C_\nu$ be the cube in $[0,1]^k$ centred at $\vv x_0$ of
sidelength $l^{-\nu}$. Then
$$
|A\cap C_\nu|=\int_{C_\nu}\chi_A(\vv x)d\vv
x\stackrel{(\ref{e:008})}{\le}\int_{C_\nu}\chi_A\Big(l^\nu\vv
x+\frac{\vv s}{l}\Big)d\vv x=\frac{1}{l^{\nu
k}}\int_{[0,1]^k}\chi_A(\vv x) d\vv x=|C_\nu|\cdot|A|\,.
$$
Since $\vv x_0$ is a density point of $A$ and $\diam C_\nu\to0$ as $\nu \to \infty$,
the left hand side  of the above equality is asymptotically
$|C_\nu|$. Therefore, $|A|=1$.

\vspace{-3ex}

\hfill $\boxtimes$

\vspace{1.5ex}

Given a ball $B=B(x,r)$ and a real number $c>0$, we denote by $cB$
the `scaled' ball  $B(x,cr)$. The next lemma  is a basic covering
result from geometric measure theory usually referred to as the
$5r$-lemma. For further details and proof the reader is refereed
to \cite{Mattila-1995}.

\begin{lemma}\label{lemma1}
Every collection $\cC$ of balls of uniformly bounded diameter in a
metric space $\Omega$ contains a disjoint subcollection $\cG$ such
that \ $$ \bigcup_{B\in\cC}B\subset\bigcup_{B\in\cG}5B\,. $$
\end{lemma}

\vspace{1.5ex}

We immediately make use of the covering lemma to show that the
Lebesgue measure of a reasonably general $\limsup$ set is
unchanged with respect to `scaling' by a constant factor.

\begin{lemma}\label{prop1}
Let $\{S_i\}_{i\in\N}$ be a sequence of subsets in $[0,1]^k$,
$\{\delta_i\}_{i\in\N}$ be a sequence of positive numbers such
that $\delta_i\to0$ as $i\to\infty$ and let
$$
\Delta(S_i,\delta_i):=\{x\in[0,1]^k:\dist(S_i,x)<\delta_i\}\,.
$$
Then for any real number $C>1$, the sets
$$A:=\limsup_{i\to\infty}\Delta(S_i,\delta_i) \qquad {\rm  and
} \qquad B:=\limsup_{i\to\infty}\Delta(S_i,C\delta_i)$$ have the
same Lebesgue measure.
\end{lemma}

\noindent\textit{Proof.} First of all notice that the sets
$\Delta(S_i,\delta_i)$ are open and therefore Lebesgue measurable.
Since $C>1$ we have that $A\subset B$. For each $i \in \N$,  let
$\cB_i$ denote the collection of balls $\{B(x,\delta_i):x\in
S_i\}$. Thus, we have that $
\Delta(S_i,\delta_i)=\bigcup_{B\in\cB_i}B\,. $ By
Lemma~\ref{lemma1}, there is a disjoint subcollection $\cG_i$ of
$\cB_i$ such that
\begin{equation}\label{a}
\bigcup^\circ_{B\in\cG_i}B\subset
\Delta(S_i,\delta_i)=\bigcup_{B\in\cB_i}B\subset
\bigcup_{B\in\cG_i}5B\,.
\end{equation}
 Since $S_i\subset[0,1]^k$ we have that
every ball $B\in\cG_i$ is contained in the cube
$[-\delta_i,1+\delta_i]^k$. It follows that $\cG_i$ is a finite
disjoint collection of balls.

\noindent If $z\in\Delta(S_i,C\delta_i)$, then there is a $y\in
S_i$ such that $|z-y|<C\delta_i$. Furthermore, by (\ref{a}) there
exists  a ball $B=B(x,\delta_i)\in\cG_i$ such that $y\in 5B$.
Therefore, $|z-x|\le|z-y|+|y-x|<(5+C)\delta_i$. Thus we have shown
that
\begin{equation}\label{b}
\Delta(S_i,C\delta_i)\subset \bigcup_{B\in\cG_i}(5+C)B\,.
\end{equation}
Now given a constant $\lambda>0$, let $
 \cC(\lambda):=\limsup_{i\to\infty}\bigcup_{B\in\cG_i}\lambda B\,.
$ This is the set of $x$ such that $x\in \lambda B$ for some
$B\in\cG_i$ for infinitely many $i$. Then, (\ref{a}) and (\ref{b})
imply that
\begin{equation}\label{c}
\cC(1)\subset A\subset B\subset \cC(5+C)\,.
\end{equation}
By Lemma~\ref{lemma0}, the sets $\cC(\lambda)$ with $\lambda>0$
have the same Lebesgue measure irrespective of $\lambda$.
Therefore, in view of  (\ref{c}) the sets $A$ and $B$ must have
the same Lebesgue measure.

\hfill $\boxtimes$

\section{Proof of Theorem \ref{t1}}

On following the arguments of \cite{Gallagher-61:MR0133297},  it
is easily verified  that $\, \anyA =  [0,1]^{nm}\,$ if the
following condition
\begin{equation}\label{e:009}
\text{$\Psi(\vv q)/|\vv q|\to0$ as $|\vv q|\to\infty$}
\end{equation}
is  violated. Therefore, without loss of generality we assume
that (\ref{e:009}) is satisfied.

When considering $\anyA$, the error of approximation is rigidly
determined by the function $\Psi$. In proving  Theorem \ref{t1},
it is extremely useful to introduce a certain degree of
flexibility within the error of approximation. Given $\anyA$, let
$$
\anyF=\bigcup_{k=1}^\infty\cA^\circ_{n,m}(k\Psi)\,.
$$
Clearly,  $\anyF \supset \anyA$. However, as a consequence of
Lemma \ref{prop1} we have that
\begin{equation} \label{sv100}
|\anyF| \, = \,
|\anyA| \ .
\end{equation}
 Clearly, Theorem~\ref{t1} follows on establishing the analogous
statement for $\anyF$.

\begin{theorem}\label{t2+}
For any $n,m$ and $\Psi$ we have that $|\anyF|\in\{0,1\}$\,.
\end{theorem}

\subsection{Proof of Theorem \ref{t2+}}

We establish the theorem by considering the sets $\FnmP$, $\FdnmP$
and $\FddnmP$ separately.

\bigskip

\noindent\textbf{The set $\FnmP$\,:} Clearly,  $\FnmP$ is
invariant under the translation $T: \vv X\mapsto 2\vv X\pmod 1$.
Thus,  the desired statement for the set  $\FnmP$ is a trivial
consequence of Lemma~\ref{l4}.
\vspace{-3ex}

\hfill $\boxtimes$

\vspace{1.5ex}

\bigskip

\noindent\textbf{The set $\FdnmP$\,: } By definition,  $\FdnmP$
consists of points $\vv X\in[0,1]^{nm}$ for which  there exists a
constant $C=C(\vv X)>0$ such that
\begin{equation}\label{e:010}
    |\vv q\vv X+\vv p|<\ C\,\Psi(\vv q)\qquad\text{and}\qquad     \gcd(\vv p,\vv q)=1
\end{equation}
holds for infinitely many $(\vv p,\vv q)\in\Z^m\times\Z^n\nz$.
Now, for each prime $l$ consider the following subsets of
$\FdnmP$:
$$
\begin{array}{l}
 \cS_0(l)=\{\vv X\in[0,1]^{nm}:\text{$\exists \  C>0$ so that (\ref{e:010}) holds for i.m.
 $(\vv p,\vv q)$ with }l\nmid d=\gcd(\vv q)\},\\[2ex]
 \cS_1(l)=\{\vv X\in[0,1]^{nm}:\text{$\exists \ C>0$ so that (\ref{e:010}) holds for i.m.
 $(\vv p,\vv q)$
 with }l\,\|\,d=\gcd(\vv q)\},\\[2ex]
\cS_2(l)=\{\vv X\in[0,1]^{nm}:\text{$\exists \ C>0$ so that
(\ref{e:010}) holds for i.m. $(\vv p,\vv q)$ with
}l^2\,|\,d=\gcd(\vv q)\} .
\end{array}
$$
Here  $l\|d$ means that $l$ divides $d$ but $l^2$ does not divide
$d$. Note that
\begin{equation}\label{e:012}
    \FdnmP=\cS_0(l)\cup\cS_1(l)\cup\cS_2(l)\,.
\end{equation}

\noindent Suppose  $\vv X\in\cS_0(l)$. Then (\ref{e:010}) is satisfied
for infinitely many $(\vv p,\vv q)$ with $l\nmid d=\gcd(\vv q)$.
On setting $\vv q':=\vv q$ and $\vv p':=l\vv p$,  we have  that
\begin{equation*}
    |\vv q'(l\vv X)+\vv p'|<lC\Psi(\vv q')
\end{equation*}
holds for infinitely many $(\vv p',\vv q')\in\Z^m\times\Z^n\nz$
with
\begin{equation}\label{e:014}
    \gcd(\vv p',\vv q')=1\,.
\end{equation}
The coprimality condition is readily verified by making use of the fact
that \mbox{$l\nmid\gcd(\vv q)$}. Thus, if  $\vv X\in\cS_0(l)$ then  $l\vv
X\in\cS_0(l)$. Therefore the set $\cS_0(l)$ is invariant under the
transformation $T: \vv X\mapsto l\vv X\pmod 1$ and  Lemma~\ref{l4} implies that
$|\cS_0(l)|$ is 0 or 1.

\noindent For $j\in\{1,\dots,n\}$, let $ \cS_{1,j}(l)$ denote the
set  of $\vv X\in[0,1]^{nm}$ such that (\ref{e:010}) is satisfied
for infinitely many $(\vv p,\vv q)$ with $l\|q_j$. Recall, that
$q_j$ is the $j$'th coordinate of $\vv q=(q_1,\dots,q_n)  $.
Clearly,
$$\cS_1(l)=\bigcup_{j=1}^n\cS_{1,j}(l)   \, . $$
\noindent Suppose  $\vv X \in \cS_{1,j}(l)$ for some
$j\in\{1,\dots,n\}$. Let $\vv S_j\in\Z^{nm}$ denote the integer
matrix with zero entries everywhere except in the $j$-th row where
every entry is $1$. Then $\vv q\vv S_j=(q_j,\dots,q_j) \in \Z^m$.
By definition,  (\ref{e:010}) is satisfied for infinitely many
$(\vv p,\vv q)$ with $l\|q_j$. On setting $\vv q':=\vv q$ and $\vv
p':=l\vv p-\frac1l\vv q\vv S_j$, we have that
\begin{equation}\label{e:015}
    \left|\vv q'\Big(l\vv X+\frac1l\vv
S_j\Big)+\vv p'\right|<lC\Psi(\vv q')\,
\end{equation}
holds for infinitely many $(\vv p',\vv q')\in\Z^m\times\Z^n\nz$
satisfying (\ref{e:014}). Thus, if $\vv X \in \cS_{1,j}(l)$  then
  $l\vv X+\frac1l\vv
S_j\in\cS_{1,j}(l)$. Therefore the set $\cS_{1,j}(l)$ is invariant
under the transformation $$T \, : \, \vv X\mapsto l\vv
X+\frac1l\vv S_j\pmod 1\,$$ and Lemma~\ref{l4} implies that
$|\cS_{1,j}(l)|$ is 0 or 1. Thus, $\cS_1(l)$ is a finite union of
sets with measure 0 or 1 and so  $|\cS_1(l)|$ is also $0$ or $1$.

In view of (\ref{e:012}), the upshot of the above results for
$|\cS_0(l)|$ and $|\cS_1(l)|$ is that if there exists a prime $l$
such that $\cS_0(l)$ or $\cS_{1}(l)$ is of positive measure then
$|\FdnmP| = 1$. Thus, without loss of generality we can assume
that such a prime does not exist and so  by (\ref{e:012}) we have
that
\begin{equation}\label{sv:101}
|\cS_2(l)|=|\FdnmP|  \quad  \text{ for every prime $l$ } \, .
\end{equation}

\noindent Suppose $\vv X\in\cS_2(l)$ and fix any $\vv
S\in\Z^{nm}$. Then (\ref{e:010}) is satisfied for infinitely many
$(\vv p,\vv q)$ with $l^2|d=\gcd(\vv q)$. On setting $\vv q':=\vv
q$ and $\vv p':=\vv p -\frac 1l \, \vv q\vv S$, we have that
\begin{equation}\label{e:016}
    \left|\vv q'\left(\vv X+\frac1l\vv S\right)+\vv p'\right|<\ C\,\Psi(\vv
    q)\,
\end{equation}
holds for infinitely many $(\vv p',\vv q')\in\Z^m\times\Z^n\nz$
satisfying (\ref{e:014}).  Thus, if $\vv X \in \cS_{2}(l)$  then
  $\vv X+\frac1l\vv S\in\cS_{2}(l)$ for any $\vv S\in\Z^{nm}$. Therefore
  the set
$\cS_{2}(l)$ is invariant under any transformation $\vv X\mapsto
\vv X+\frac1l\vv S  \pmod 1 $ with $\vv S\in\Z^{nm}$. In other
words, $\cS_2(l)$ is $\frac1l$-periodic in every coordinate. Thus,
for any cube $\cC_l$ in $[0,1]^{nm}$ of sidelength $\frac1l$ we
have that
$$
     |\cS_2(l)\cap \cC_l|=|\cS_2(l)| \cdot |\cC_l|\,.
$$
In view of (\ref{sv:101}),  it follows that
\begin{equation}\label{e:017}
     |\FdnmP\cap \cC_l|=|\FdnmP| \cdot |\cC_l|\,
\end{equation}
for any prime $l$. Now suppose that $|\FdnmP|>0$. Then there is a
density point $\vv X_0$ of $\FdnmP$. For each prime $l$, let
$\cC_l$ denote the  cube in $[0,1]^{nm}$ centred at $\vv X_0$ of
sidelength $\frac1l$. Then,
$$
|\FdnmP\cap \cC_l|\sim|\cC_l| \qquad\text{as } \ \ l\to\infty\,.
$$
This together with (\ref{e:017}) implies that $|\FdnmP|=1$ and
thereby completes the proof of  Theorem~\ref{t1} for the set
$\FdnmP$. \vspace{-3ex}

\hfill $\boxtimes$

\bigskip

\noindent\textbf{The set $\FddnmP$\,: } To establish the desired
zero-one statement for the set $\FddnmP$,  we modify in the
obvious manner the argument given above for the set $\FdnmP$.
Naturally, ``$\gcd(p_j,\vv q)=1$ for all $j=1,\dots,m$''  will
replace ``$\gcd(\vv p,\vv q)=1$'' appearing in  (\ref{e:010}).
Similarly, the condition that ``$\gcd(p'_j,\vv q')=1$ for all
$j=1,\dots,m$'' will replace the coprimality condition
(\ref{e:014}). The rest  remains pretty much unchanged.
\vspace{-3ex}

\hfill $\boxtimes$


\section{Further results and questions}

\noindent\textbf{$\Psi$-well approximable points.} The
various sets  of \emph{$\Psi$-well approximable points} are
defined by requiring that the constant $C>0$ appearing in
(\ref{e:010}) can be made arbitrarily small.
 More precisely,
$$
\anyW:=\bigcap_{k=1}^\infty\cA^\circ_{n,m}(k^{-1}\Psi)\,.
$$

Lemma~\ref{prop1} readily implies the following statement.

\begin{theorem}\label{t3}
For any $n,m$ and $\Psi$ we have that $|\anyW|=|\anyA|$.
\end{theorem}

Theorem \ref{t3} combined with (\ref{sv100}) and Theorem~\ref{t2+}
trivially implies the zero-one law for $\Psi$-well approximable
sets.

\begin{corollary}
For any $n,m$ and $\Psi$ we have that $|\anyW|\in\{0,1\}$\,.
\end{corollary}

\medskip

\noindent\textbf{$\Psi$-badly approximable points.} Naturally, the various
sets of \emph{$\Psi$-badly
approximable} points can be thought of as being complementary to the $\Psi$-well
approximable sets.    More precisely,
$$
\anyB:=\anyA\setminus\anyW\,.
$$

An immediate consequence of  Theorem~\ref{t3} is the following
result.

\begin{corollary}\label{svcor}
For any $n,m$ and $\Psi$ we have that $|\anyB|=0$\,.
\end{corollary}

\noindent The classical set of badly approximable real numbers
$\mathbf{Bad}:=\cB_{1,1}(q\mapsto q^{-1})$ is known to have full
Hausdorff dimension; i.e.  $\dim\textbf{Bad}=1$. For a general
function $\psi:\N\to\Rp$ with various mild growth conditions,
Bugeaud \cite{Bugeaud-03:MR2006007}, answering a question posed in
\cite{Beresnevich-Dickinson-Velani-01:MR1866488}, has shown that
$\cB_{1,1}(\psi)$ has full Hausdorff dimension; i.e.
$\dim\cB_{1,1}(\psi)=\dim\cA_{1,1}(\psi) $.  It view of this and
Corollary \ref{svcor} it is reasonable to ask the following
question.

\begin{question}\label{q1}
Is it true that the set $\BnmP$ has full Hausdorff dimension; i.e.
$$\dim\BnmP=\dim\AnmP \, ? $$
\end{question}
A weaker form of this question in which  $\BnmP$ is replaced by
$\AnmP\setminus\Anm(\Psi')$ with $\Psi'(\vv q)=o(\Psi(\vv q))$ as
$|\vv q|\to\infty$, can be found in
\cite{Beresnevich-Dickinson-Velani-01:MR1866488}.  Note that if
answer to the above question is yes, then we automatically have
that  $ \dim\anyB=\dim\anyA $.

\bigskip

\noindent\textbf{Multi-error approximation.} Observe that the
inequality given by  (\ref{e:001}) can be rewritten as a system of
$m$ inequalities; namely
$$|\vv q\vv X^{(j)}+p_j|<\Psi(\vv q) \ \quad j=1,\dots,m\,,
$$ where $\vv X^{(j)}$ is the $j$-th column of $\vv X$. Thus,
the error of approximation associated with each linear form is
determined by $\Psi$ and is the same. More generally, we consider
the system
\begin{equation}\label{e:018}
  |\vv q\vv X^{(j)}+p_j|<\Psi_j(\vv q) \ \qquad j=1,\dots,m\, ,
\end{equation}
with $\Psi_j:\Z^n\to\Rp$ and so the error of approximation is
allowed to differ from one linear form to the next. Let
$\Anm^\circ(\Psi_1,\dots,\Psi_m)$ denote the  `multi-error'
analogue of $\anyA$  -- obtained by replacing  (\ref{e:001}) with
(\ref{e:018}) in the definition of  $\anyA$ . Naturally, this
enables us to define the  multi-error analogues of  $\anyF$,
$\anyW$ and $\anyB$.

\noindent Without much effort, it is possible to establish the
multi-error analogue of Theorem~\ref{t2+} -- the proof is
practically unchanged.

\begin{theorem}
For any $n,m$ and $\Psi_1,\dots,\Psi_m$ we have that $|\Fnm^\circ(\Psi_1,\dots,\Psi_m)| \, \in\{0,1\}$.
\end{theorem}

\noindent If the statement of
Lemma~\ref{prop1}  can be  generalized to the multi-error framework the above theorem  would
answer the following question and thereby yield the
analogue of Theorem~\ref{t1}.
\begin{question}
Is it true that
$|\Fnm^\circ(\Psi_1,\dots,\Psi_m)| =|\Wnm^\circ(\Psi_1,\dots,\Psi_m)|$?
\end{question}

\noindent Note that if the answer to Question~2 is yes, then so is
the answer to our next question.

\begin{question}
Is it true that $|\Bnm^\circ(\Psi_1,\dots,\Psi_m)|=0$?
\end{question}

\bigskip

\noindent\textbf{Multiplicative approximation.} Given
$\Psi:\Z^n\to[0,+\infty)$, let $\cA^\times_{n,m}(\Psi)$ be the set
of $\vv X\in[0,1]^{nm}$ such that
\begin{equation}\label{e:019}
  \prod_{i=1}^m\|\vv q\vv X^{(j)}\|<\Psi(\vv q) \,
\end{equation}
holds for infinitely many $\vv q\in\Z^n$. Here  $\|\cdot\|$
denotes the distance to the nearest integer.  Naturally, this
enables us to define the  associated multiplicative
 sets $\cF^\times_{n,m}(\Psi)$,
$\cW^\times_{n,m}(\Psi)$ (multiplicatively $\Psi$-well
approximable points) and $\cB^\times_{n,m}(\Psi)$
(multiplicatively $\Psi$-badly approximable points). Clearly, if
$\Psi:=\Psi_1\cdots\Psi_m$ then
$$\cA(\Psi_1,\dots,\Psi_m)\subset\cA^\times_{n,m}(\Psi)\, ,$$
$$\cF(\Psi_1,\dots,\Psi_m)\subset\cF^\times_{n,m}(\Psi)\, ,$$
$$\cW(\Psi_1,\dots,\Psi_m)\subset\cW^\times_{n,m}(\Psi) \, . $$
However, it is easily seen that
$$\cB(\Psi_1,\dots,\Psi_m)\not\subset\cB^\times_{n,m}(\Psi) \, .$$

\begin{question}
Is it true that $\cA^\times_{n,m}(\Psi)$, $\cF^\times_{n,m}(\Psi)$ and $\cW^\times_{n,m}(\Psi)$
are of measure $0$ or $1$?
\end{question}
\begin{question}
Is it true that $|\cB^\times_{n,m}(\Psi)|=0$?
\end{question}
Note that when $n=1$, $m=2$ and  $\Psi(q):=q^{-1}$, the answer to Question 5  is  for obvious reasons yes. Indeed,
it is conjectured that
$$\cB^\times_{1,2}(q\mapsto q^{-1})=\emptyset \ . $$
This  is   Littlewood's famous
conjecture in the theory of Diophantine approximation.

\medskip

\noindent\textbf{Approximation by rational planes.} The inequality
given by (\ref{e:001}) takes on two `extreme' forms of rational
approximation depending on whether $n=1 $ or $m=1$. When $m=1$, it
corresponds to  approximating arbitrary points by
$(n-1)$-dimensional rational planes  (i.e. rational hyperplanes)
and gives rise to the dual theory of Diophantine approximation.
When $n=1$, it corresponds to approximating  arbitrary points by
$0$-dimensional rational planes (i.e. rational points) and gives
rise to the simultaneous theory of Diophantine approximation. For
$d \in \{0, \ldots n-1\}$, it is natural to consider the
Diophantine approximation theory in which points in $\R^n$ are
approximated by $d$-dimensional rational planes -- the dual and
simultaneous theories just represent the extreme. The foundations
have been  developed in some depth by W.M.~Schmidt
\cite{Schmidt-67:MR0213301} in the sixties and more recently by
M.~Laurent \cite{Laurent}. However,  apart from the extreme cases,
there appears to be no analogue of Theorem~\ref{t1} within the
theory of approximation by $d$-dimensional rational planes.

\bigskip

\noindent\textbf{Approximation by algebraic numbers.} Sprindzuk's
\cite{Sprindzuk-1979-Metrical-theory} celebrated proof of Mahler's
conjecture \cite{Mahler-1932b} led to Baker
\cite{Baker-1966-On-a-theorem-of-Sprindzuk} making the following
stronger conjecture that was eventually established by Bernik
\cite{Bernik-1989}. \emph{Let $n\in\N$ and $\psi:\N\to(0,+\infty)$
be monotonic. Then for almost every real $x$ the inequality
\begin{equation}\label{e:021}
    |P(x)|<H(P)^{-n+1}\psi(H(P))
\end{equation}
holds for finitely many $P\in\Z[x]$ with $\deg P\le n$ if }
\begin{equation}\label{e:022}
    \sum_{r=1}^\infty\psi(r) <\infty\,.
\end{equation}
Here $H(P)$ is the height of $P$; i.e.  the maximum of the
absolute values of the coefficients of $P$. The case
$\psi(h):=h^{-1-\ve}$ corresponds to Mahler's conjecture. In
\cite{Beresnevich-99:MR1709049} it has been shown that if the sum
in (\ref{e:022}) diverges and $\psi$ is monotonic, then for almost
every real $x$ inequality (\ref{e:021}) holds infinitely often.
More recently \cite{Beresnevich-05:MR2110504}, the monotonicity
assumption in Bernik's convergence result has been removed.
However, removing the monotonicity assumption  from the divergence
result remains an open problem  akin to the Duffin-Schaeffer
conjecture. In the first instance,  it would be natural and
desirable to
 ask for a  zero-one law.
\begin{question}
Is it true that the set of $x\in[0,1]$ such that $(\ref{e:021})$
holds for infinitely many $P\in\Z[x]$ with $\deg P\le n$ is of
measure $0$ or $1$?
\end{question}
The following is a related question concerning explicit
approximation by algebraic numbers.
\begin{question}
Is it true that the set of real $x\in[0,1]$ such that
\begin{equation}\label{e:023}
    |x-\alpha|< H(\alpha)^{-n}\psi(H(\alpha))
\end{equation}
holds for infinitely many real algebraic $\alpha$ of
$\deg\alpha\le n$ is of measure $0$ or $1$?
\end{question}
Here $H(\alpha)$ stands for the height of the minimal defining
polynomial of $\alpha$. If (\ref{e:022}) is satisfied, a simple
application of the Borel-Cantelli lemma shows that the set under
consideration is of measure zero. On the other hand, if the sum in
(\ref{e:022}) diverges and $\psi$ is monotonic the set under
consideration is known to have measure one -- see
\cite{Beresnevich-99:MR1709049}. The upshot is that  Question~7 only
needs to be considered when $\psi$ is non-monotonic and the sum in
(\ref{e:022}) diverges.

\vspace{6ex}

\noindent{\em Acknowledgements. } We would like to thank Wolfgang
Schmidt for his many wonderful theorems that have greatly influenced
our research. Also, thank you for your generous support over the
years. Finally and most importantly -- happy number seventy five
Wolfgang!

\bigskip


{\small

\def\cprime{$'$} \def\cprime{$'$} \def\cprime{$'$} \def\cprime{$'$}
  \def\cprime{$'$} \def\cprime{$'$} \def\cprime{$'$} \def\cprime{$'$}
  \def\cprime{$'$} \def\cprime{$'$} \def\cprime{$'$} \def\cprime{$'$}
  \def\cprime{$'$} \def\cprime{$'$} \def\cprime{$'$} \def\cprime{$'$}
  \def\cprime{$'$}

}

{\small
\vspace{5mm}

\noindent Victor V. Beresnevich: Department of Mathematics,
University of York,

\vspace{-2mm}

\noindent\phantom{Victor V. Beresnevich: }Heslington, York, YO10
5DD, England.


\noindent\phantom{Victor V. Beresnevich: }e-mail: vb8@york.ac.uk

\vspace{5mm}

\noindent Sanju L. Velani: Department of Mathematics, University of York,

\vspace{-2mm}

\noindent\phantom{Sanju L. Velani: }Heslington, York, YO10 5DD, England.


\noindent\phantom{Sanju L. Velani: }e-mail: slv3@york.ac.uk

}

\end{document}